\documentclass[10pt,twocolumn,twoside]{IEEEtran}

\usepackage{setspace}

\pdfoutput 1

\usepackage{graphicx}          
\usepackage{amsmath,amssymb,amsthm} 
\usepackage[english]{babel}
\usepackage{amsfonts}
\usepackage{amstext}
\usepackage{caption}
\usepackage{tabularx}
\usepackage{savesym}
\savesymbol{AND}
\usepackage{algorithm}
\usepackage{algorithmic}
\usepackage{color}
\usepackage{todonotes}
\usepackage{subfig}

\usepackage{exStyle}

\newtheorem{theorem}{Theorem}[section]
\newtheorem{lemma}[theorem]{Lemma}

\newtheorem{definition}[theorem]{Definition}
\newtheorem{assumption}[theorem]{Assumption}

\newtheorem{problem}[theorem]{Problem}
\newcommand{\pp}[1]{\ensuremath{\mathbb P\left\{#1\right\}}}
\newcommand{\ve}[1]{{\bf #1}}

\newcommand{\C}{\mathcal C}
\newcommand{\Pp}{\mathcal P}
\newcommand{\E}{\mathcal E}
\newcommand{\V}{\mathcal V}

\begin{document}
%

\author{Roel Dobbe*, Ye Pu, Jingge Zhu, Kannan Ramchandran, Claire Tomlin
\thanks{All authors are with the Department of Electrical Engineering and Computer Sciences at UC Berkeley. *Corresponding author: dobbe@berkeley.edu.}}

\title{Customized Local Differential Privacy for Multi-Agent Distributed Optimization}

\maketitle



%


\begin{abstract}
Real-time data-driven optimization and control problems over networks may require sensitive information of participating users to calculate solutions and decision variables, such as in traffic or energy systems. 
Adversaries with access to coordination signals may potentially decode information on individual users and put user privacy at risk.
We develop \emph{local differential privacy}, which is a strong notion that guarantees user privacy regardless of any auxiliary information an adversary may have, for a larger family of convex distributed optimization problems. 
The mechanism allows agent to customize their own privacy level based on local needs and parameter sensitivities.
We propose a general sampling based approach for determining sensitivity and derive analytical bounds for specific quadratic problems.
We analyze inherent trade-offs between privacy and suboptimality and propose allocation schemes to divide the maximum allowable noise, a \emph{privacy budget}, among all participating agents.
Our algorithm is implemented to enable privacy in distributed optimal power flow for electric grids.

\end{abstract}

\section{Introduction}

Advances in sensing and computing enable various infrastructures, such as traffic or energy networks, to perform optimization and control problems in real-time throughout a network.
Often the scale of such problems desires a distributed implementation that can be solved quickly enough to allow for high frequency control actions.
To enable this, a network may be split up into sub-networks governed by different agents, who exchange their local optimization variables with neighbors and/or a central operator to iteratively solve the optimization problem.
Exchanging optimization variables between agents and the changes therein may reveal private information, such as whether someone is home and what kind of appliances someone is using \cite{hart_nonintrusive_1992}. 
In addition, there is growing understanding that secondary information may be inferred from the communicated variables, including the parameters used in the local objective and constraints, which may reveal sensitive information such as prices and capacity~\cite{han_differentially_2017}.

To make matters more challenging, different agents may be competing with each another to serve an operator with their service.
Knowing the control capacity of and prices negotiated by other players can help in negotiating with the operator and leads to strategic behavior and untruthful communication, which harms the quality of solution to the distributed optimization problem.
As such, both personal privacy and commercial privacy needs require the development of agent-to-agent distributed optimization algorithms that can mask sensitive information in objectives and constraints. 

In recent years, various privacy-preserving algorithms have been proposed for distributed optimization and control problems, using various privacy metrics.
The differential privacy framework \cite{dwork_algorithmic_2014} has gained most attention, and is particularly lauded for its robustness to auxiliary side information that an adversary might have to complement information gained from a particular algorithm, providing stronger privacy guarantees than other existing metrics. 
The framework assumes a setting in which sensitive information is stored in a database by a trustworthy curator, which can provide answers to external queries.
A system is made differentially private by randomizing its answers in such a way that the distribution over published outputs is not too sensitive to changes in the stored data. 
Perturbation can be designed to make it provably difficult for an adversary to make inferences about individual records from the published outputs.

In the setting of distributed optimization, each agent is its own curator managing its own locally private information and communication of its optimization variables to neighboring agents or a central operator. 
In order to preserve differential privacy, each curator has to ensure that the output of queries, that is the communicated variables, remain approximately unchanged if local parameters relating to its objective or constraints are modified.

\subsection*{Related Work}
This work complements an existing and rapidly growing body of literature on incorporating \emph{differential privacy} into resource allocation and, most relevant here, in distributed optimization, control and networked systems.
Earlier work by Hsu et al.~\cite{hsu_privately_2014} developed differential privacy-preserving algorithms for convex optimization problems that are solved in a \emph{central} fashion, considering general linear programs (LPs) with either private objectives or constraint. Dong et al.~\cite{dong_differential_2015} consider privacy in a game theoretic environment, motivated by traffic routing in which the origins and destinations of drivers are considered private. Jia et al.~\cite{jia_privacy-enhanced_2017} consider occupancy-based HVAC control a thend treat the control objective and the location traces of individual occupants as private variables, using an information-theoretic privacy metric.
A recent elaborate tutorial paper by Cort\'{e}s et al.~\cite{cortes_differential_2016} covers differential privacy for \emph{distributed} optimization, and distinguishes between objective-perturbing and message-perturbing strategies for distributed optimization.
In the first category, each agent’s objective function is perturbed with noise in a differentially private manner, which guarantees differential privacy at the functional level and is preferred for systems with asymptotically stable dynamics~\cite{nozari_differentially_2018}.
In the second category, coordination messages are perturbed with noise before sent, either to neighbors or a central node, depending on the specific algorithm. Huang et al.~\cite{huang_differentially_2015} proposed a technique for disguising private information in the local objective function for distributed optimization problems with strongly convex separable objective functions and convex constraints. Han et al.~\cite{han_differentially_2017} considered problems where the private information is encoded in the individual constraints, the objective functions need to be convex and Lipschitz continuously differentiable, and the constraints have to be convex and separable.
Other related works are Mo and Murray~\cite{mo_privacy_2017} who aim to preserve privacy of agents' initial states in average consensus and Katewa et al.~\cite{katewa_privacy_2017} who explore an alternative trade-off between privacy and the value of cooperation (rather than performance) in distributed optimization.  In \cite{shoukry_privacy-aware_2016}, a privacy-aware optimization algorithm is analyzed using the cryptography notion of zero knowledge proofs.  More recently, \cite{zhu_differentially_2018} considers differentially private algorithms over time-varying directed networks.



The above works are selective in that these consider privacy-preserving mechanisms for constraints, objectives or initial states. 
An exception is the work Hsu et al.~\cite{hsu_privately_2014} on LPs, which can handle both private objectives and constraints.
We complement the literature by proposing a mechanism that preserves private objectives and constraints for optimization problems with strongly convex objectives and convex constraints.

\subsection*{Contributions}
Motivated by personal and commercial privacy concerns in distributed OPF, we investigate the problem of \emph{preserving differential privacy of local objectives and constraints in distributed constrained optimization with agent-to-agent communication}. Compared to previous works on privacy-aware distributed optimization~\cite{han_differentially_2017,huang_differentially_2015}, we consider the notion of  \textit{local differential privacy}, which is a refined (and more stringent) version of differential privacy~\cite{duchi_local_2013}. It allows each agent in a network to \emph{customize} its own privacy level, based on individual preferences and characteristics. Furthermore, most previous works consider privacy protection for either individual objective function~\cite{huang_differentially_2015} or the individual constraint~\cite{han_differentially_2017},  our more general formulation enables us to provide privacy guarantees on both local objective function parameters and local constraint parameters. Specifically, the proposed algorithm solves a general class of convex optimization problems where each agent has a local objective function and a local constraint, and agents communicate with neighbors/adjacent agents with no need for a central authority.

We show that the private optimization algorithm can be formulated as an instance of the Inexact Alternating Minimization Algorithm (IAMA) for Distributed Optimization \cite{pu_inexact_2014}.
This algorithm allows provable convergence under computation and communication errors. 
This property is exploited to provide privacy by injecting noise large enough to hide sensitive information, while small enough to exploit the convergence properties of the IAMA.
We derive the trade-off between the privacy level and sub-optimality of the algorithm.
The trade-off between sub-optimality and differential privacy allows us to determine a \emph{privacy budget} that captures the allowable cumulative variance of noise injected throughout the network that achieves a desired level of (sub-)optimality. We propose two pricing schemes to ensure fair and efficient allocation of the privacy budget over all participating/bidding DER owners.,

\section{Preliminaries and Problem Statement}

In this section, we consider a distributed optimization problem on a network of $M$ sub-systems (nodes). The sub-systems communicate according to a fixed undirected graph $G =(\mathcal{V},\mathcal{E})$. The vertex set $\mathcal{V} = \{1,2,\cdots,M\}$ represents the sub-systems and the edge set $\mathcal{E}\subseteq \mathcal{V}\times \mathcal{V}$ specifies pairs of sub-systems that can communicate. If $(i,j)\in \mathcal{E}$, we say that sub-systems $i$ and $j$ are neighbors, and we denote by $\mathcal{N}_i = \{j| (i,j)\in \mathcal{E}\}$ the set of the neighbors of sub-system $i$. Note that $\mathcal{N}_i$ includes $i$. The cardinality of $\mathcal{N}_i$ is denoted by $|\mathcal{N}_i|$. We use a vector $[v]_i$ to denote the local variable of subsystem $i$ and $[v]_i$ can be of different dimensions for different $i$. 
The collection of these local variables is denoted as $v=[v^{T}_1,\cdots,v^{T}_M]^{T}$. 
Furthermore, the concatenation of the local variable $[v]_i$ of sub-system $i$ and the variables of its neighbors $[v]_j, j\in\mathcal N_i$ is denoted by $z_{i}$. 
With appropriate selection matrices $E_i$ and $F_{ji}$, the variables have the following relationship: $z_{i} = E_iv$ and $[v]_i = F_{ji}z_{j}$, $j\in\mathcal{N}_i $, which implies the relation between the local variable $[v]_i$ and the global variable $v$, i.e. $[v]_i = F_{ji}E_jv$, $j\in\mathcal{N}_i $. We consider the following distributed optimization problem:
%
\begin{problem}[Distributed Optimization]
\begin{align}
\min_{z,v}  & \quad \sum^{M}_{i=1} f_{i}(z_i)\\
s.t. &\quad  z_i\in \mathbb{C}_{i}, \quad z_i = E_iv, \quad i=1,2,\cdots ,M \,,
\end{align}
\label{pr:dist_opt}
\end{problem}
\noindent where $f_i$ is the local cost function for node $i$ which is assumed to be strongly convex  with a convexity modulus $\rho_{f_i}>0$, and to have a Lipschitz continuous gradient with a Lipschitz constant $L(\nabla f_i)>0$.   The constraint $\mathbb{C}_i$ is assumed to be a convex set which represents a convex local constraint on $z_{i}$, i.e. the concatenation of the variables of sub-system $i$ and the variables of its neighbors. 

The above problem formulation is fairly general and can represent a large class of problems in practice. In particular  it includes the  following quadratic programming problem, which we study as a particular instance in our applications.
\begin{problem}[Distributed Quadratic Problem]
\begin{align}
\min_{z,v}  & \quad \sum^{M}_{i=1} z^{T}_iH_iz_i + h^{T}_iz_i \label{pr:quadratic}\\
s.t. &\quad  C_iz_i\leq c_i, \quad z_i = E_iv, \quad i=1,2,\cdots ,M\,, \nonumber
\end{align}
\label{pr:dist_quad_opt}
\end{problem}
\noindent where $H_i \succ 0$.  In particular, we will assume that the smallest eigenvalue of $H_i$ satisfies $\lambda_{min}(H_i) := \lambda_{min}^{(i)} > 0$.

\subsection{Local Differential Privacy}


We present definitions and properties for differential privacy. 
Let $\Pp, \Pp'$ be two databases with private parameters in some space $\mathcal X$ containing information relevant in executing an algorithm. 
Let $\text{adj}:\mathcal X\times\mathcal X\mapsto [0,\infty)$ denote a metric defined on $\mathcal X$ that encodes the adjacency or distance of two elements. 
A mechanism or algorithm $\mathcal A$ is a mapping from $\mathcal X$ to some set denoting its output space.


\begin{definition}[Differential Privacy]. A randomized algorithm $\mathcal{A}$ is $\epsilon$-differential private if for all $\mathcal{S}\subseteq \mbox{range}(\mathcal{A})$ and for all database $\Pp, \Pp'$ satisfying $\text{adj}(\Pp,\Pp')\leq 1$, it holds that
\begin{equation}
\mbox{Pr}[\mathcal{A}(\Pp) \in \mathcal{S}] \leq e^{\epsilon}\cdot\mbox{Pr}[\mathcal{A}(\Pp') \in \mathcal{S}] \,,
\end{equation}
where the probability space is over the mechanism $\mathcal{A}$. 
\label{def: Differential privacy}
\end{definition}

\noindent This definition of differential privacy is suitable for cases where one uniform level of privacy needs to across all elements in the databases. 
We now consider a distributed algorithm $\mathcal A(\Pp_1,\ldots, \Pp_M)$ in a network with $M$ agents for solving an optimization problem in a collaborative way, where $\Pp_i$ denotes the private parameters of agent $i$. The outputs of the mechanism are the message exchanged between nodes in the network over the time horizon of iterations. 
This mechanism induces $M$ local mechanisms $\mathcal A_1(\Pp_1,\ldots, \Pp_M),\ldots, \mathcal A_M(\Pp_1,\ldots, \Pp_M)$, each executed by one agent. The output of one local mechanism $\mathcal A_i$ is the message sent out by node $i$, i.e. $\text{range}(\mathcal A_i)\subseteq\text{range}(\mathcal A)$. It is important to realize that although one local mechanism, say $\mathcal A_i$, does not necessarily have direct access to the input/database $\Pp_j, j\neq i$ of other nodes, the output of $\mathcal A_i$ could still be affected by $\Pp_j, j\neq i$ because of the interactions among different nodes. For this reason, we explicitly write $\Pp_1,\ldots, \Pp_M$ as input to all local mechanisms.

We now let each agent~$i$ specify its own level of privacy~$\epsilon_i$.
To formalize this specification, we require a definition:
\begin{definition}[Local Differential Privacy]
Consider a (global) mechanism $\mathcal A$ for a network with $M$ nodes, and  $M$ local mechanisms $\mathcal A_i, i=1,\ldots, M$ induced by $\mathcal A$.  We say that the mechanism $\mathcal A$ is \textit{$\epsilon_i$-differentially locally private for node $i$}, if for any $\mathcal S_i\in\text{range}(\mathcal A_i)$ it satisfies that
\begin{align}
\frac{\pp{\mathcal A_i(\Pp_1,\ldots, \Pp_i,\ldots, \Pp_M)\in\mathcal S_i}}{ \pp{\mathcal A_i(\Pp_1,\ldots, \Pp_i',\ldots, \Pp_M)\in\mathcal S_i}}\leq e^{\epsilon_i} \,,
\end{align}
where $\text{adj}(\Pp_i,\Pp_i')\leq 1$. Moreover, we say that the mechanism $\mathcal A$ is \textit{$(\epsilon_1,\ldots, \epsilon_M)$-differentially private}, if  $\mathcal A$ is $\epsilon_i$-differentially locally private for all nodes, where $i=1,\ldots, M$.
\label{def:local_privacy}
\end{definition}

\begin{figure}[!hbt]
\centering
\includegraphics[width = 0.45\textwidth]{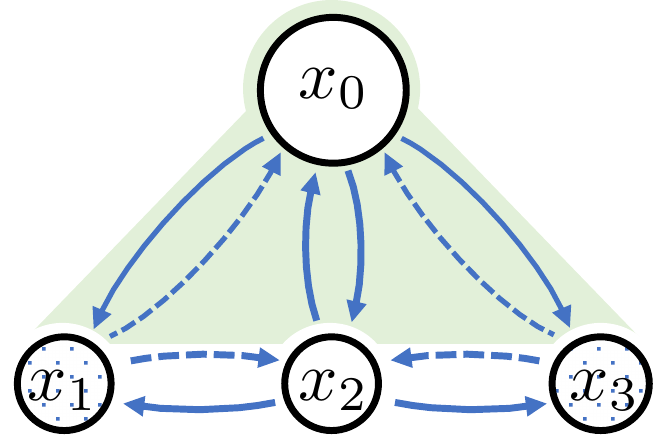}
\caption{A system with 4 subsystems. Arrows denote message directions in distributed optimization. Here, node 0 denotes a central node that communicates with all other nodes, which can be left out (shaded part). The other 3 nodes represent subsystems. Nodes that have a privacy specification are patterned with dots and send out messaged perturbed by noise, as indicated by dashed arrows.}
\label{fig:concept}
\end{figure}

Figure~\ref{fig:concept} presents the concept of local differential privacy pictorially, showing the various considerations that can be taken when designing for local/customized privacy.
Firstly, one may desire to include a central node 0 that communicates with all subsystems or implement a fully distributed problem between the subsystems that does not rely on any central node. 
The former will lead to better convergence properties as information spreads more easily throughout the network, the latter will benefit privacy by making it harder to collect information from across the network.
Regardless, the method allows for the privacy to be purely local, strengthening the notion developed in~\cite{han_differentially_2017}, which implements differential privacy in distributed optimization through a trusted central node, assuming a star-shaped communication structure with no agent-to-agent communication.
Secondly, subsystems may have varying levels of privacy. In  Figure~\ref{fig:concept}, subsystem 1 and 3 have a local privacy specification, while subsystems 0 and 2 do not. The systems with local differential privacy have outgoing messages perturbed by noise as indicated by the dashed arrows.
As such, the method is flexible to various forms of distributed optimization or control problems with heterogeneous privacy and control properties across its nodes/agents, extending the work in \cite{huang_differentially_2015}, which considers a similar fully distributed algorithm but specifies a uniform privacy level for all agents.

\section{Main Results}
\label{sec:main_results}
\subsection{Differentially Private Distributed Optimization}
\label{subsec:dpdo}

In this section, we describe a distributed optimization algorithm  (Algorithm \ref{al: Distributed optimization with errors}) for solving Problem~\ref{pr:dist_opt} with privacy guarantee, based on the results in \cite{pu_inexact_2014}. To solve the optimization problem in a distributed way,  we split the overall problem into small  local problems according to the physical couplings of the sub-systems. This approach gives the resulting algorithm the desired feature that each node only needs to communicate with its neighbors and the computations can be performed in parallel for every subsystem. To guarantee local differential privacy, a noise term is added to the message at each time, before it is sent out to other nodes.


\begin{scriptsize}
\begin{algorithm}
\caption{Differentially private distributed algorithm}
\begin{algorithmic} 
\REQUIRE Initialize $\mu^{0}_i =0 \in \mathbb{R}^{z_i}$, $\tau^{0} = \min_{1\leq i \leq M} \{\rho_{f_i}\}$ and $\tau^{k} = \frac{1}{\tau^{0}k}$
\FOR {$k=1,2,\cdots$}
\STATE 1: $z^{k}_i = \mbox{argmin}_{z_i\in \mathbb{C}_{i}}\lbrace  f_i(z_i) + \langle \mu^{k-1}_{i}, -z_i\rangle\rbrace + \delta^{k}_i$
\STATE 2: Send $z^{k}_i$ to all the neighbors of agent $i$.
\STATE 3: $[v^{k}]_{i} = \frac{1}{|\mathcal{N}_i|} \sum_{j\in\mathcal{N}_i} [z^{k}_j]_i$.
\STATE 4: Send $[v^{k}]_{i}$ to all the neighbors of agent $i$.
\STATE 5: $\mu^{k}_i = \mu^{k-1}_i + \tau^{k}(E_i v^{k}-z^{k}_i)$
\ENDFOR
\end{algorithmic}
\label{al: Distributed optimization with errors}
\end{algorithm}
\end{scriptsize}

We start with defining the private parameters of agent~$i$ to be the collection of parameters for its local objective function and constraints in Problem~\ref{pr:dist_opt}, $\mathcal P_i:=(f_i, \mathbb C_i)$. 
Algorithm~\ref{al: Distributed optimization with errors} can be seen as a \textit{global mechanism} $\mathcal A$ which takes input data $\mathcal P_1,\ldots,\mathcal P_M$ and produces output $z_j^k, [v^k]_j$ for all $j=1,\hdots, M$ and $k=1,2,\hdots$.   In particular, the output of $\mathcal A$ up to iteration $K$ is given by
\begin{align}
\mathcal A^{K}(\mathcal P_1,\ldots,\mathcal P_M):=(\ve z^1,\ldots, \ve z^K, [\ve v]^1, \ldots, [\ve v]^K) \,,
\end{align}
where we use the bold letter $\ve z^k$ to denote the collection of  $(z_1^k,\ldots,z_M^k)$ at time $k$ and $[\ve v]^k$ to denote the collection of $([v]_1^k\ldots, [v]_M^k)$ at time $k$.
Recall that we are mainly concerned with local privacy at different nodes, we use $\mathcal A_i$ to denote the \textit{local mechanism} induced by $\mathcal A$ at node $i$, which takes input data $\mathcal P_1,\ldots, \mathcal P_M$ and produces output $z_i^k, [v^k]_i$ for $k=1,2\ldots$, namely
\begin{align}
\mathcal A_i^K(\mathcal P_1,\ldots,\mathcal P_M):=( z_i^1,\ldots, z_i^K, [v]_i^1, \ldots, [v]_i^K) \,,
\label{eq:Ai_output}
\end{align}
which injects noise~$\delta_i^k$ to the solution $z_i^k$ of its local problem, according to some distribution.

The definition of the distance $d(\mathcal P_i,\mathcal P_i')$ between to two problems $\mathcal P_i$ and $\mathcal P_i'$ depends on specific applications. For the quadratic problem in (\ref{pr:quadratic}), the local objective $f_i$ and local constraint $\mathbb C_i$ are parametrized by the matrix and vectors $H_i, h_i, C_i, c_i$. One possible definition of the distance for this special case is given by weighted sum of matrix/vector norms
\begin{align}
\text{adj}(\mathcal P_i,\mathcal P_i')&=a_1\|H_i-H_i'\|+a_2\|h_i-h_i'\| \nonumber\\
&+a_3\|C_i-C_i'\|+a_4\|c_i-c_i'\|  \,,
\label{eq:distance_quadratic}
\end{align}
for some weights $a_1,a_2,a_3,a_4$. The choice of the weights will put more emphasis sensitivity on a specific matrix or vector and its corresponding parameters. 
These can be used to define the ``region of adjacency'', which represents the parameter differences that are considered to be meaningful for protecting with a differential privacy mechanism. The choice of the norm can also change the behavior, for instance by treating all parameters equal (such as for the $\ell_1$-norm or the Frobenius norm) or only worrying about the maximum distance (such as for the $\ell_{\infty}$-norm). For some applications, only certain entries of the matrices and vectors represent private information. In this case, the distance should be defined only with respect to these ``private entries''. 

In the classical setup of differential privacy where the database representation is considered, there is a natural definition of adjacency: namely two databases are adjacent if they differ by only one element. When extending the adjacency definition to the space of functions, there does not exist a natural candidate for adjacency. Even with a choice such as (\ref{eq:distance_quadratic}), it is not completely clear how to normalize the norms. We should point out that this is a common situation encountered in similar problem setups such as  \cite{han_differentially_2017,huang_differentially_2015}. However, resolving this issue would require an explicit connection between the privacy level $\epsilon$ and a concrete specification supplied by the application (such as in~\cite[Section B]{cortes_differential_2016}), which is out of scope for this paper.



\subsection{Privacy analysis}


In this section, we derive the differential privacy level with Algorithm~\ref{al: Distributed optimization with errors}.   To state the main result, we first define  
\begin{align}
\displaystyle
g(\mathcal P_j,\mu) :=\arg \min_{z\in \mathbb{C}_{j}}\lbrace  f_j(z) + \langle \mu, -z\rangle\rbrace \,,
\label{eq:function_g}
\end{align}
where $\mathcal P_j$ encapsulates the local objective function $f_j$ and constraints $\mathbb C_j$. The main result is given in the following theorem.
\begin{theorem}[Local Differential Privacy]
Consider Algorithm \ref{al: Distributed optimization with errors} for solving Problem \ref{pr:dist_opt}. Assume that each element of the noise vector $\delta_i^k$ in Algorithm \ref{al: Distributed optimization with errors} is independently chosen according  the Laplace distribution with the density function $p(\delta_i^k)=\frac{1}{2\sigma_i^k}\exp(-\|\delta_i^k\|/\sigma_i^k)$ for $i=1,\ldots,M$ and $k=1,\ldots, K$. 
Then for \textit{all} $i=1,\ldots, M$, this algorithm is locally $\epsilon_i$-differentially private for node $i$ where 
\begin{align}
\epsilon_i=\Theta_i \sum_{k=1}^K \frac{1}{\sigma_i^k} \,,
\label{eq:privacy_level}
\end{align}
and 
\begin{align}
\Theta_i:=\displaystyle \max_{d(\mathcal P_i,\mathcal P_i')\leq 1, \mu} \|g(\mathcal P_i,\mu)-g(\mathcal P_i',\mu)\|
\label{eq:sensitivity}
\end{align}
is called the \textit{sensitivity} of the optimization problem.
\label{thm:privacy}
\end{theorem}
%
\noindent 
\begin{IEEEproof}
The proof is given in the Supplementary Material.
\end{IEEEproof}

\subsection{Sensitivity Calculation}
\label{sec:sensitivity}
In order to evaluate the privacy level $\epsilon_i$ provided in Theorem \ref{thm:privacy},  we need to calculate the sensitivity $\Theta_i$.  This calculation is itself an optimization problem which can be written as follows.
\begin{align}
\max  \ & \Theta_i:= \|z^*-z'^{*}\| \,,\label{prob:calculate_sensitivity}\\
\text{s.t.} \ & z^* = \displaystyle \arg \min_{z\in \mathbb C_i}{f_i(z)+\langle \mu, -z\rangle} \,, \nonumber \\
& z'^{*} = \displaystyle \arg \min_{z\in \mathbb C_i'}{f_i'(z) + \langle \mu, -z\rangle} \,, \nonumber \\
& \text{adj}(\mathcal P_i(f_i,\mathbb C_i), \mathcal P_i'(f_i',\mathbb C_i')) \leq 1 \,, \nonumber \\
&\text{with variables }(z^*,z'^{*},f_i,f_i',\mathbb C_i,\mathbb C_i', \mu) \,. \nonumber
\end{align}
This problem belongs to the class of bi-level optimization problems, for which there is in general no efficient algorithm to find the global optimal solution. In the rest of this section, we will specialize our problem to the the class of quadratic problems in (\ref{pr:quadratic}) and provide some refined analysis.  In this case, we represent the optimality condition of $z^*,z'^{*}$ in terms of the KKT condition of the optimization problem $\min_{z\in\mathbb C_i} f_i(z)+\langle \mu, -z\rangle$. Hence, the  optimization problem (\ref{prob:calculate_sensitivity}) can be rewritten explicitly
\begin{align}
\max \ & \Theta_i:= \|z^*-z'^{*}\| \,, \label{eq:sensitivity_quadratic}\\
\text{s.t.} \ &  z^*\in \text{KKT}(H_i,h_i,C_i,c_i,w_i, \mu) \,, \label{cons1} \\
& z'^{*}\in\text{KKT}(H_i',h_i',C_i',c_i',w_i', \mu) \,, \label{cons2} \\
& a_1\|H_i-H_i'\|+a_2\|h_i-h_i'\|\nonumber\\
& +a_3\|C_i-C_i'\|+a_4\|c_i-c_i'\| \leq 1 \,, \label{cons3} \\
&\text{with variables } (z^*,z^{'*}, H_i, H_i', h_i, h_i', \nonumber\\
&C_i, C_i', c_i, c_i', w_i, w_i', \mu ) \,, \nonumber
\end{align}
where the set $\text{KKT}(H_i,h_i,C_i,c_i,w_i,\mu)$ is defined as
\begin{align*}
\text{KKT}(H_i,h_i,C_i,c_i,w_i,\mu):=\{z | &  H_iz+h_i-\mu+C_i^Tw_i = 0 \,, \\
& C_iz_i\leq c_i \,, \\
&w_i\geq 0 \,, \\
&w_{i,j}(C_iz-c_i)_j = 0 \text{ for all }j\} \,,
\end{align*}
where $w_i$ represents Lagrangian multipliers for the optimization problem (\ref{prob:calculate_sensitivity}). 
Notice that we used the distance $d(\Pp_i,\Pp_i')$ as defined for the quadratic problem in~\eqref{eq:distance_quadratic}.

We state our first observation regarding the sensitivity of quadratic problems as formulated in Problem~\ref{pr:dist_quad_opt}.
\begin{lemma}[Sensitivity for Problem~\ref{pr:dist_quad_opt}]
If we specialize Problem  \ref{pr:dist_opt} to the problem in (\ref{pr:quadratic}) with the distance $d(\mathcal P_i,\mathcal P_i')$ defined in (\ref{eq:distance_quadratic}), $\Theta_i$ can be simplified as
\begin{align}
\Theta_i:=\max_{\text{adj}(\mathcal P_i,\mathcal P_i')\leq 1} \|g(\mathcal P_i,\mu=0)-g(\mathcal P_i',\mu=0)\| \,,
\label{eq:quad_case_sensitivity}
\end{align}
that is we lose explicit dependency on lagrange variable $\mu$.
\label{lemma:independent_mu}
\end{lemma}

In the rest of this section, we discuss how to estimate the sensitivity for  generic quadratic problems using a sample-based method. We point out that a more general discussion of this approach can be found in \cite{calafiore_scenario_2006}. For the sake of notation, we write $\Theta_i$ in (\ref{eq:sensitivity_quadratic}) as  $\Theta_i(\Pp)$ where $\Pp$ denotes the collection of all variables $(z^*,z'^{*}, H_i, H_i', \tilde h_i, \tilde h_i',C_i, C_i', c_i, c_i', w_i, w_i')$. 
Furthermore we use $\mathbb C_{\Theta}$ to denote the polynomial constraints given by (\ref{cons1}), (\ref{cons2}) and (\ref{cons3}). With these notations, the  optimization problem in (\ref{eq:sensitivity_quadratic}) can be rewritten as
\begin{align}
\min \ & \gamma \,, \label{prob:equivalent}\\
\text{s.t.} \ & \Theta(\Pp) - \gamma\leq 0, \  \forall \Pp \in \mathbb C_{\Theta} \,. \nonumber
\end{align}
The idea of sample-based approach is to randomly draw many instances of $\Pp$ from the set $\mathbb C_{\Theta}$, and find the maximum $\Theta(\Pp)$ using these samples.  Namely, we solve the following problem
\begin{align}
\gamma^N:=\min \ & \gamma \,, \label{prob:estimate}\\
\text{s.t.} \ & \Theta(\Pp_s) - \gamma\leq 0, \ \forall \Pp_s, \ s=1,\ldots, N \,, \nonumber
\end{align}
where $\Pp_1,\ldots, \Pp_N$ are $N$ randomly drawn samples from $\mathbb C_{\Theta}$. More specifically for our problem, we randomly sample the parameters $ H_i, H_i', \tilde h_i, \tilde h_i',C_i, C_i', c_i, c_i', w_i, w_i'$ from their sets. For each sampled set of parameters, we solve the original optimization problem (\ref{pr:quadratic}) to obtain $z^*$ and $z'^{*}$, hence obtain one estimate  $\Theta(\Pp_s):=\|z^*-z'^{*}\|$. After $N$ samples, the maximal $\Theta(\Pp_s)$ is set to be $\gamma^N$, which gives a lower bound on the sensitivity. To quantify the quality of this approximation, the following definition is introduced.
\begin{definition}[Random Sampling, Definition 1, 2 in~\cite{calafiore_uncertain_2005}]
Let $\gamma^*$ denote the optimal solution to problem (\ref{prob:equivalent}) and $\gamma^N$ be a candidate solution retrieved through solving~\eqref{prob:estimate}. 
We say $\gamma^N$ is an $\alpha$-level robustly feasible solution, if $V(\gamma^N)\leq \alpha$, where $V(\gamma^N)$ is defined as
\begin{align}
V(\gamma^N):=\mathbb P\{ \Pp \in \mathbb C_{\Theta}:  \gamma^N \le \gamma^*\} \,.
\end{align}
where the measure is induced from the uniform sampling of the constraint set.
\end{definition}
In other words, $V(\gamma^N)$ is the portion of $\Pp$ in $\mathbb C_{\Theta}$ that was not explored after $N$ samples, which, if explored, would yield a higher function value than $\gamma^N$ and thereby a tighter lower bound on sensitivity~$\Theta(\Pp)$. 
The following result relates the number of samples $N$ to the quality of the approximation.

\begin{lemma}[Sampling Rule, Corollary 1 in~\cite{calafiore_uncertain_2005}]
\label{le:sampling_rule}
For a given $\alpha \in [0,1]$ and $\beta \in [0,1]$ and let $N\geq \frac{1}{\alpha\beta}-1$. Then with probability no smaller than $1-\beta$, the solution $\gamma^N$ given by (\ref{prob:estimate}) is a solution with $\alpha$-level robust feasibility for Problem (\ref{prob:equivalent}).
\end{lemma}

\noindent This result gives the minimum number of samples to find an approximate optimal solution to~\eqref{prob:estimate} with high probability (larger than $1-\beta$).

\subsection{Analytical Upperbound for Sensitivity in Special Cases}

\noindent By making additional assumptions on the problem statement for the Distributed Quadratic Problem, formulated in Problem~\ref{pr:quadratic}, we derive an analytical upper bound on the sensitivity $\Theta_i$.
\begin{lemma}[Special Case - Protecting $H_i$ in Problem~\ref{pr:dist_quad_opt}]
\label{le:sensitivity_H}
Assume $f_i(z_i):=\frac{1}{2}z_i^TH_ix+h_i^Tz_i$ and the distance between two problems is defined as $\text{adj}(\mathcal P_i,\mathcal P_i')=\|H_i-H_i'\|_2$. That is, the privacy requirement only concerns about the matrix $H_i$.  Also assume that the local variable $z_i$ is bounded as $\|z_i^k\|\leq G_i$, then
\begin{equation}
\Theta_i\leq \frac{G_i}{\lambda_{\text{min}}^{(i)}} \,,
\label{eq:sensitivity_Hi}
\end{equation}
where $\lambda_{\text{min}}^{(i)}$ is the lower bound on the eigenvalues  defined in Problem~\ref{pr:dist_quad_opt}.
\label{lemma:sensitivity_change_quadratic}
\end{lemma}

%
%

\begin{lemma}[Special Case - Protecting $h_i$ in Problem~\ref{pr:dist_quad_opt}]
\label{le:sensitivity_h}
Assume $f_i(z_i):=\frac{1}{2}z_i^TH_ix+h_i^Tz_i$ and the distance between two problems is defined as $\text{adj}(\mathcal P_i,\mathcal P_i')=\|h_i-h_i'\|_2$. That is, the privacy requirement only concerns about the vector $h_i$. Then
\begin{align*}
\Theta_i\leq \frac{1}{\lambda_{min}^{(i)}} \,,
\end{align*}
where $\lambda_{min}^{(i)}$ is the lower bound on the eigenvalues  defined in Problem~\ref{pr:dist_quad_opt}.
\label{lemma:sensitivity_change_linear}
\end{lemma}

These lemmas give closed-form expressions of the upper bounds on $\Theta_i$ for two special cases, and they will be useful for our applications in Section \ref{sec:d-OPF}. Notice however that the upperbounds do not serve as a straighforward design principle for scaling eigenvalues to lower sensitivity. As a result of scaling the eigenvalues, the distance metric~$d(\cdot,\cdot)$ will also change; higher eigenvalues generally result in smaller sensitivity but also a narrower and less meaningful ``region of adjacency'', as discussed in Section~\ref{subsec:dpdo}, Equation~\eqref{eq:distance_quadratic}.

\subsection{Convergence properties of the distributed optimization algorithm}\label{se: DOA}

Privacy comes with a price. The noise term in  Algorithm~\ref{al: Distributed optimization with errors} makes the convergence rate  slower than the case without noise.  However, it is possible to prove that even with random noise, the algorithm converges in expectation. To show this, we first write out the dual problem of Problem \ref{pr:dist_opt} as follows.
\begin{problem}[Dual Problem of Problem \ref{pr:dist_opt}]\label{pr:dual problem of distributed optimization problem}
\begin{equation*}
\min \quad -D(w) = \underbrace{\sum^{M}_{i=1} f^{\star}_{i}(w_i)}_{\phi(w)}  + \underbrace{\mbox{I}_{\{E^{T}w=0\}}(w)}_{\psi(w)} \,,
\end{equation*}
\end{problem}
\noindent with the matrix $E := [E^{T}_1,E^{T}_2,\cdots ,E^{T}_M]^{T}$  and dual variable $w:=[w_1^T,\ldots, w_M^T]^T$. 
We use  $f_i^*$ to denote the conjugate function of $f_i: \mathbb C_i\rightarrow \mathbb R$, defined as $f_i^{\star}(w) = \sup_{z\in \mathbb C_i} (w^{T}z - f(z))$ and   $\mbox{I}_{\mathbb{S}}$ denotes the indicator function on a set $\mathbb{S}$
\begin{equation*}
\mbox{I}_{\mathbb{S}}(x) =  \begin{cases} 0 &\mbox{if } x \in \mathbb{S} \\ 
\infty & \mbox{if } x \notin \mathbb{S} \end{cases}.
\end{equation*}

\noindent The dual problem is of the form:
\begin{problem}[Dual Problem Form]
\label{pr:problem ISTA}
\begin{align*}
\min_{w\in \mathbb W} \quad \Phi(w) =\phi(w)+\psi(w) \enspace .
\end{align*}
\end{problem}

\noindent The stochastic proximal-gradient method (stochastic PGM), as given in Algorithm \ref{al:stochastic PGM}, is a method to solve the dual problem above.  
\begin{scriptsize}
\begin{algorithm}
\caption{Stochastic Proximal-Gradient Method}
\begin{algorithmic}
\REQUIRE Require $w^{0} \in \mathbb W$ and step size $\tau^{k}<\frac{1}{\rho_{\phi}k}$
\FOR {$k=1,2,\cdots$}
\STATE 1: $w^{k} =  \mbox{prox}_{\tau^{k} \psi}(w^{k-1} - \tau^{k} (\nabla \phi(w^{k-1}) + e^k))$
\ENDFOR
\end{algorithmic}
\label{al:stochastic PGM}
\end{algorithm}
\end{scriptsize}
The proximity operator in Algorithm \ref{al:stochastic PGM} is defined as
\begin{align*}
\mbox{prox}_{\tau \psi}(v):=\text{argmin}_{x} \left( \psi(x) + (1/2\tau\|x-v\|^2)\right) \,.
\end{align*}
This algorithm has been studied extensively~\cite{rakhlin_making_2012,shamir_stochastic_2013} and requires the following assumptions.
\begin{assumption}[Assumptions for Problem~\ref{pr:problem ISTA}] \text{  } \label{as:ISTA with strong convexity}
\begin{itemize}
\item  $\phi$ is a strongly convex function with a convexity modulus $\rho_{\phi}>0$, and has a Lipschitz continuous gradient with a Lipschitz constant $L(\nabla \phi)>0$.\color{black}
\item $\psi$ is a lower semi-continuous convex function, not necessarily smooth.
\item  The norm of the gradient of the function $\phi$ is bounded, i.e., $\|\nabla \phi(w)\|^{2}\leq B^{2}$ for all $w\in\mathbb{W}.$
\item The variance of the noise $e^{k}$ is equal to $ \sigma^{2}$, i.e., $\mathbb{E}[\|e^{k}\|^{2}]\leq \sigma^{2}$ for all $k$.
\end{itemize}
\end{assumption}



The key observation in our proof of convergence is the following lemma, showing that executing Algorithm \ref{al: Distributed optimization with errors} on the original problem \ref{pr:dist_opt} is equivalent to executing Algorithm \ref{al:stochastic PGM} on the dual problem \ref{pr:dual problem of distributed optimization problem}.
\begin{lemma}[Equivalence Dual Problem]\label{le:the equivalence between distributed optimization algorithm and stochastic PGM}
Consider using Agorithm \ref{al: Distributed optimization with errors} on  Problem \ref{pr:dist_opt} and using Algorithm \ref{al:stochastic PGM} on Problem \ref{pr:dual problem of distributed optimization problem}. Further assume that  Algorithm~\ref{al: Distributed optimization with errors} is initialized with the sequence $\mu_j^0, z_j^0$, for $j=1,\ldots, M $, and Algorithm~\ref{al:stochastic PGM} is initialized with the sequence $w_j^0$ where $w_j^0=\mu_j^0$ for $j=1,\ldots, M$.  Then $w_j^k=\mu_j^k$ for all $k=1,2\ldots$ and all $j=1,\ldots, M$, and the error terms in Algorithm~\ref{al: Distributed optimization with errors} and Algorithm~\ref{al:stochastic PGM} have the relationship $e^{k} = \delta^{k}= [\delta^{kT}_1,\delta^{kT}_2,\cdots ,\delta^{kT}_M]^{T}$.
\end{lemma}
A proof sketch is provided in the supplementary material, largely building on results in~\cite{pu_inexact_2016}.
Based on the equivalence shown in Lemma~\ref{le:the equivalence between distributed optimization algorithm and stochastic PGM}, we are ready to provide the following theorem showing the convergence properties of Algorithm~\ref{al: Distributed optimization with errors}.

\begin{theorem}[Suboptimality]
\label{th:convergence rate of inexact AMA}
Consider Algorithm \ref{al: Distributed optimization with errors}. Assume that the local variables $z_i^k$ are bounded as $\|z_i^{k}\|^{2}\leq G^{2}_i$ for all $k\geq 0$ and for all $i=1,\cdots,M$.  We have that for any $k>1$, the expected suboptimality is bounded as
\begin{equation}\label{eq: upper-bound of the distributed optimization algorithm}
\mathcal S_k := \mathbb{E}[|D(w^{k}) - D(w^{\star})|]\leq \frac{4\sum^{M}_{i=1} (G^{2}_i+\sigma^{2}_i)}{\rho^{2}_{\phi}k}\enspace .
\end{equation}
\end{theorem}
\begin{IEEEproof}
Consider applying Algorithm \ref{al:stochastic PGM} on Problem \ref{pr:problem ISTA} with Assumption \ref{as:ISTA with strong convexity}, we can first show that
\begin{equation}
\mathbb{E}[\|w^{k} - w^{\star}\|]\leq \frac{4(B^{2}+\sigma^{2})}{\rho^{2}_{\phi}k}\enspace .
\end{equation}
The proof of this claim follows the same flow as the proof of Theorem~1 in \cite{rakhlin_making_2012} by noticing the following two facts: if the function $f$ is convex, closed and proper, then 
\begin{equation}
 \|\mbox{prox}_{f}(x) - \mbox{prox}_{f}(y)\| \leq \|x-y\| \enspace ,
\end{equation}
and the variance of the gradient of $\phi$ is bounded by $\mathbb{E}[\|\nabla \phi(w^{k-1}) + e^k\|^{2}]\leq B^2 + \sigma^{2}$ for all $k>1$.

We then apply this result to our Problem \ref{pr:dual problem of distributed optimization problem}.  The gradient of the first objective is equal to:
\begin{align*}
\nabla_{w} \phi(w^{k}) &= \nabla_{w} \sum^{M}_{i=1} f^{\star}_{i}(w^{k}_i) = [z^{k^{T}}_{1}, \cdots , z^{k^{T}}_{M}]^{T}
\end{align*}
with $w^{k} = [w^{k^{T}}_{1}, \cdots , w^{k^{T}}_{M}]^{T}$. With our assumption in the theorem, the dual gradient is bounded as:
\begin{equation*}
\|\nabla \phi(w^{k})\|^{2} \leq \sum^{M}_{i=1} \|z^{k}_{i}\|^{2} \leq \sum^{M}_{i=1} G^{2}_i \enspace .
\end{equation*}
Finally the claim follows because of the equivalence result in Lemma~\ref{le:the equivalence between distributed optimization algorithm and stochastic PGM}.
\end{IEEEproof}

Theorem \ref{eq: upper-bound of the distributed optimization algorithm} shows that the suboptimality gap for Algorithm~\ref{al: Distributed optimization with errors} is bounded above by a function that is linear in the number of agents $M$ and linear in the noise variances~$\sigma_i\,,i = 1,\hdots,M$. The convergence rate is of order $\mathcal O\left(\dfrac{1}{k}\right)$. This is an improvement compared to the motivating work by Han et. al~\cite{han_differentially_2017}, which achieved a convergence rate of $\mathcal O\left(\dfrac{1}{\sqrt{k}}\right)$, and is similar to the algorithm proposed by Huang et. al~\cite{huang_differentially_2015}.

\section{Application: Distributed Optimal Power Flow}\label{sec:d-OPF}
This section presents a simplified optimal power flow (OPF) problem that inspires the proposed control approach. We consider the setting of a radial distribution feeder, and consider the flow of real power on its branches. We formulate the power flow model and the OPF objectives and develop the distributed OPF problem according to the quadratic problem, as defined in~\eqref{pr:quadratic}. We then discuss the parameters that are subject to privacy requirements and interpret the trade-offs developed in Section~\ref{sec:main_results}.

\subsection{Simplified Optimal Power Flow}
Solving the simplified OPF problem requires a model of the electric grid describing both topology and impedances. This information is represented as a graph $\mathcal G = (\V,\E)$, with $\V$ denoting the set of all buses (nodes) in the network, and $\E$ the set of all branches (edges).
For ease of presentation and without loss of generality, here we introduce part of the linearized power flow equations over \emph{radial} networks, also known as the \emph{LinDistFlow} equations \cite{baran_optimal_1989}. In such a network topology, each bus $j$ has one upstream parent bus $\{i \ | \ (i,j) \in \E\}$ and potentially multiple downstream child buses $\{k \ | \ (j,k) \in \E\}$. By~$\mathcal D_j$ we denote the set of all buses downstream of branch~$(i,j)$.
We assume losses in the network to be negligible and model the power flowing on a branch as the sum of the downstream net load:
\begin{equation}
P_{ij} \approx \sum_{k \in \mathcal D_j} \{ p_{k}^{\text{c}} - p_{k}^{\text{g}} + u_k \}
\end{equation}
In this model,  capital $P_{ij}$ represents real power flow on a branch from node $i$ to node $j$ for all branches $(i,j) \in \E$, lower case $p_i^c$ is the real power consumption at node $i$, and $p_i^g$ is its real power generation. 
This nodal consumption and generation is assumed to be uncontrollable.
In addition, we consider controllable nodal injection $u_i$, available at a subset of nodes $i \in \C \subset \V$ that have a Distributed Energy Resource (DER).
In this case study, we aim to prevent overload of real power flow over certain critical branches in an electric network.
This aim is formulated through constraints
\begin{equation}
\begin{array}{rcl}
\displaystyle \sum_{k \in \mathcal D_j} \{ p_{k}^{\text{c}} - p_{k}^{\text{g}} + u_k \} - \overline{P_{ij}} &\le& 0 \,,\\
\displaystyle \underline{P_{ij}} - \sum_{k \in \mathcal D_j} \{ p_{k}^\text{c} - p_{k}^\text{g} + u_k \}  &\le& 0 \,, \forall (i,j) \in \mathcal E_{\text{safe}} \,,
\end{array}
\label{eq:safety_constraints}
\end{equation}
$\mathcal E_{\text{safe}} \subset \mathcal E$ denotes a subset of branches for which power flow limitations are defined, $\overline{P_{ij}},\underline{P_{ij}}$ denoting the upper and lower power flow bounds
on branch $(i,j) \in \mathcal E_{\text{safe}}$.
In addition, each controlled node $i$ is ultimately limited by the local capacity on total apparent power capacity,
\begin{equation}
 \underline{u}_i \le u_i \le \overline{u}_i \,, \, \forall i \in \C \,.
\label{eq:InvCapd}
\end{equation}
We consider a scenario in which the operator negotiates different prices for different capacities, potentially at different points in time, with different third party DER owners.
Let $u_i$ refer to the real power used for the optimization scheme from agent $i$, and $\pi_i$ denotes the price for procuring a kWatt from agent $i$. 
The optimal power flow determines the control setpoints that minimizes an economic objective subject to operational constraints. 
\begin{IEEEeqnarray}{Rl}
\min_{u_i \,, i \in \C} \hspace{8pt} & \sum_{i \in \mathcal{C}} \pi_i (u_i)^2 \,,
\IEEEyesnumber \label{eq:LinOPF} \\
\text{s.t.} \hspace{8pt} & \eqref{eq:safety_constraints} \,, \eqref{eq:InvCapd} \,.  \nonumber
\end{IEEEeqnarray}
The OPF problem \eqref{eq:LinOPF} can be recast as an instance of the quadratic distributed optimization problem~\eqref{pr:quadratic}. 
First, note that the objective is quadratic in the optimization variables $u_i$, and separable per node. 
Second, for all nodes $i \in \V$, the capacity box constraints \eqref{eq:InvCapd} are linear and fully local. 
The safety constraints \eqref{eq:safety_constraints} require communication to and computation by a central trusted node.  
To ensure strong convexity of the local problems, the economic cost objectives are shared between each agent $i$ and the central trusted node. Hence, $\forall i \in \C_{\backslash \{0\}}$, the objective reads
\begin{equation}
f_i (u_i) = \frac{\pi_i}{2} (u_i)^2 \,, 
\end{equation}
with capacity constraint~\eqref{eq:InvCapd}.
The central node 0 has objective function
\begin{equation}
\begin{array}{rcl}
f_0 (z_0) &=& \displaystyle \sum_{i \in \mathcal C} \frac{\pi_i}{2} (u_i)^2 \,, \\
\end{array}
\end{equation}
with safety constraints~\ref{eq:safety_constraints}.
As such, this distributed problem assumes a star-shaped communication structure, in which the a centrally trusted node receives all $u_i, p^{\text{c}}_i , p^{\text{g}}_i$ from the agents. The agents retrieve iterates of $u_i$ from the central node and compute a simple problem with only economic cost and a local capacity constraint.

\subsection{Private Information in Distributed OPF}
We consider assigning privacy requirements to two sets of parameters; the prices $\pi_i$ that the DSO charges to different agents in the network, and the capacities $\underline{u}_i,\overline{u}_i$ available to all agents $i \in \mathcal{C}$.
Together, these parameters provide important strategic insight into the commercial position of each agent. An operator may charge different prices for different levels of commitment or for the varying value that the operator gets from the actions of a specific agent at specific time periods or places in the network. In a natural commercial context, the operator may have an interest to hide the prices to other agents. 
In addition, in a negotiation setting, a strategic agent may want to find out the capacity available by other agents in the network to adjust its bid to the operator, so as to be the first or only agent to be considered, which could lead to asymmetric and potentially unfair bidding situations. As such, in order to give all agents with capacity a fair chance to participate, there is value in hiding the capacity (and price) parameters.

To formulate this as an instance of local differential privacy, we need to define the adjacency metric for all considered parameters. In the case of both prices and capacity, this is achieved by considering the maximum range in which these parameters are expected to lie. The distance metric proposed is the $\ell_1$-norm. Given this metric, we need to define a proper adjacency relation, which determines the maximum change in a single parameter that we aim to hide with the differentially private algorithm.
\begin{definition} \textbf{(Adjacency Relation for Distributed OPF)}: For any parameter set $\Pp = \{f_i(\pi_i), \mathbb{C}_i(\overline{u}_i, \underline{u}_i) \}$ and $\Pp' = \{ f_i'(\lambda_i'), \mathbb{C}_i'(\overline{u}_i', \underline{u}_i') \}$, we have $\text{adj}(\Pp,\Pp') \le 1$ if and only if there exists $i \in [M]$ such that
\begin{equation}
\left| \pi_i - \pi_i' \right| \le \delta \pi \,, \quad \left| \overline{u}_i - \overline{u}_i' \right| \le \delta \overline{u} \,, \quad \left| \underline{u}_i - \underline{u}_i' \right| \le \delta \underline{u} \,,
\end{equation}
and $\pi_j = \pi_j'$, $\overline{u}_j = \overline{u}_j'$, $\underline{u}_j = \underline{u}_j'$ for all $j \neq i$.
\end{definition}
By setting $\delta \pi$, $\delta \overline{u}$ and $\delta \underline{u}$ respectively as the maximum price offered per unit of energy (i.e. $\bar{\pi}$ if $\pi_i \in [0, \bar{\pi}]$) and the maximum capacity in the network (i.e. $\arg \max_{i \in \mathcal{C}} \overline{u}_i$), we ensure that all parameters in the network are properly covered by the definition.

\subsection{Interpreting Trade-offs}
\label{sec:tradeoffs}
\begin{figure*}[h]
	\centering
    \includegraphics[trim = 3cm 0 3cm 0,width = \textwidth]{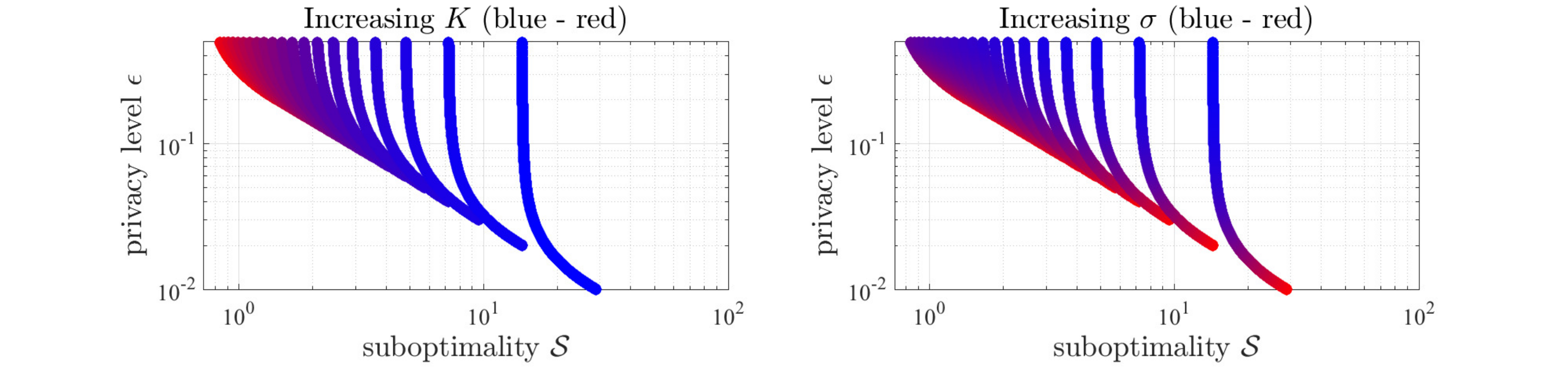}
    \caption{Achievable tradeoffs between privacy level $\epsilon$ and suboptimality $\mathcal S$, with Pareto front. (left) indicates increasing number of iterations and (right) increasing noise variance. We assume $\sigma_i = \sigma, G_i = G, \forall i \in \C$. }
    \label{fig:pareto_tradeoffs}
\end{figure*}

We analyze and interpret the theoretical results that illuminate an inherent trade-off between privacy level and suboptimality.
Assuming a fixed noise variance across all iterations, Equations~\ref{eq:privacy_level} and~\ref{eq: upper-bound of the distributed optimization algorithm}, we have the following trade-off:
\begin{equation}
\epsilon_i=\Theta_i  \frac{K}{\sigma_i}
\,, \,
\mathcal S_K \leq \frac{4\sum^{M}_{i=1} (G^{2}_i+\sigma^{2}_i)}{\rho^{2}_{\phi}K}
\end{equation}
Remember that better privacy relates to a lower privacy level~$\epsilon_i$ and a more optimal solution relates to lower suboptimality.
Unsurprisingly, an increasing number of iterations leads to worse privacy and better suboptimality. 
Conversely, a higher noise variance leads better privacy and worse suboptimality.
Figure~\ref{fig:pareto_tradeoffs} shows the region of attainable $(\epsilon,\mathcal S)$ values for the simplified OPF problem. 
The left figure shows that, for a fixed reasonable level of privacy ($\epsilon \le 10^0$), the sub-optimality will decrease for a larger number of iterations~$K$. 
The right figure shows that for a fixed level of privacy, the sub-optimality bound tightens for higher variance levels.
For a fixed level of sub-optimality, a higher noise variance~$\sigma$ achieves a lower (and hence better) privacy level.
Ideally, parameters $(\{\sigma_i\}_{i \in \C}, K)$ are chosen along the Pareto front of this graph. 
As a result, a system designer may want to define specifications,
\begin{equation}
\epsilon_i \le \overline{\epsilon}
\,, \,
\mathcal S \le \overline{\mathcal S} \,.
\end{equation}
Based on the specifications, we then want to determine feasible values for the number of iterations~$K$ and noise variance~$\sigma_1$. 
For the sake of analysis, we let all 
upper bounds be the same in the second equation, $G_i = G , \ \forall i \in \C$.
In addition, we consider the normalized noise-to-signal ratio~$\nu_i := \frac{\sigma_i}{G}$, which is more intuitive as a tunable parameter. With these steps, we can write
\begin{equation}
 \frac{K}{\nu_i} \le \frac{G_i}{\Theta_i }\overline{\epsilon}
\,, \,
\frac{M+\sum_{i=1}^M\nu_i^2}{K} \le \frac{\overline{\mathcal S}}{4} \left(\frac{ \rho_{\phi}}{G}\right)^{2}
\end{equation}
The first inequality is a local condition on the number of iterations and the (normalized) noise parameter that needs to satisfied to ensure the specified level of differential privacy for each agent. The second inequality is a collective condition that ensures the specified suboptimality.
\begin{figure*}[h]
	\centering
    \includegraphics[trim = 3cm 0 3cm 0,width = \textwidth]{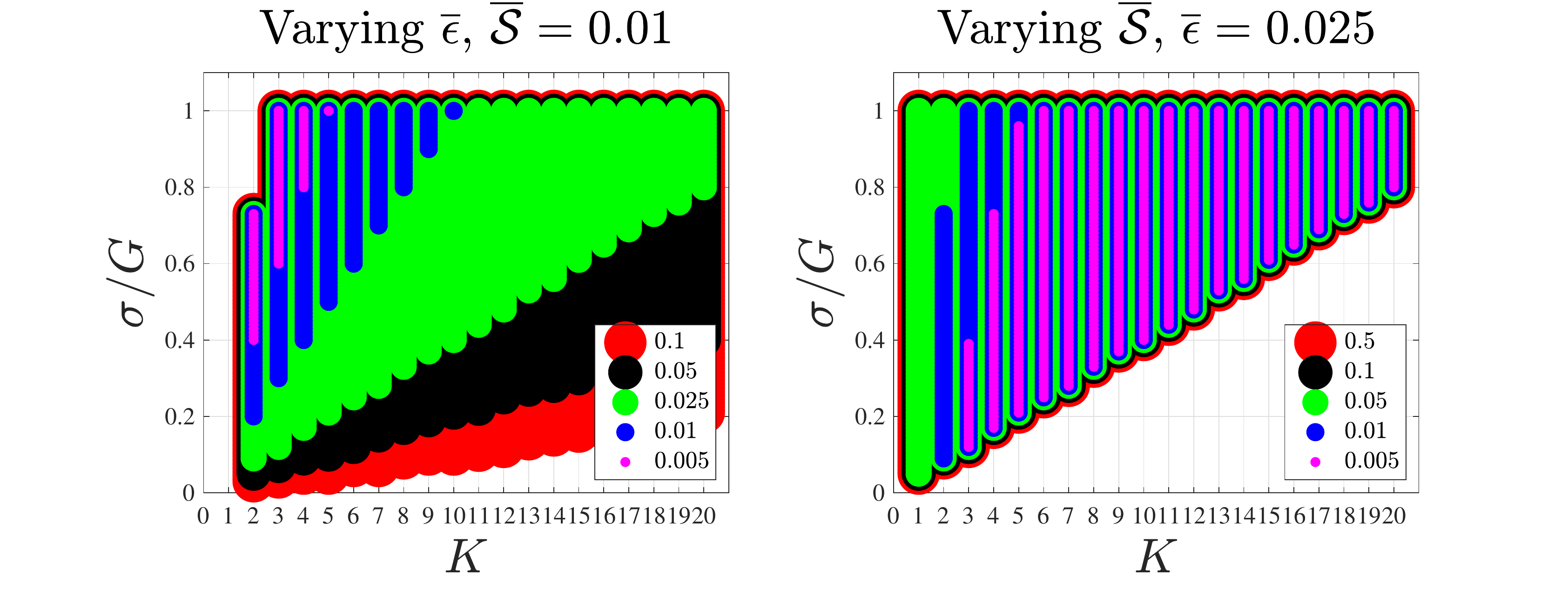}
    \caption{Feasible parameter sets $(\nu,K)$ for varying levels of $\overline{\epsilon}$ (left) and $\overline{\mathcal S}$ (right), assuming $\sigma_i = \sigma, G_i = G, \forall i \in \C$. }
    \label{fig:feasible_parameters}
\end{figure*}
These equations tell us that the suboptimality is mostly governed by $K$, as typically $0 \le \nu_i^2 \ll 1$ (in other words, we need $\sigma_i$ to be on the order of $G$ or larger to affect suboptimality).
These equations provide a specific test to determine feasible parameters~$(\nu,K)$ that satisfy a set of specifications~$(\overline{\epsilon}, \overline{\mathcal S}, M, G , \rho_{\phi}  )$. 
Note that this set may be empty if the specification are too stringent.
For sake of presentation, if we set $\sigma_i = \sigma ,\ \nu_i = \nu \ \forall i$, we get
\begin{equation}
 \frac{K}{\nu} \le \frac{G}{\Theta_i } \overline{\epsilon}
\,, \,
\frac{1+\nu^2}{K} \le \frac{\overline{\mathcal S}}{4M} \left(\frac{\rho_{\phi}}{G}\right)^{2} \ .
\end{equation}
For these settings, Figure~\ref{fig:feasible_parameters} shows the feasible set for varying levels of specifications~$(\overline{\epsilon},\overline{\mathcal{S}})$. 

We now further specify the tradeoff relations for the simplified OPF problem. 
Note that $\rho_{\phi} = \max_i \pi_i := \pi_{\text{max}}$, and $G_i = \max(|\underline{u}_i|,|\overline{u}_i|) = u_{i,\text{max}}$ and $G = \max_i G_i = u_{\text{max}}$. This yields
\begin{equation}
\frac{K}{\nu_i} \le \frac{u_{i,\text{max}}}{\Theta_i} \overline{\epsilon}_i
\ , \
\frac{M+\sum_{i=1}^M \nu_i^2}{K} \le \frac{\overline{\mathcal S}}{4} \left(\frac{\pi_{\text{max}}}{ u_{\text{max}}} \right)^{2} \,,
\end{equation}
where we have maintained the assumption that~$G_i = G , \ \forall i \in \C$ for the second inequality.
The first equation shows that the ratio of the number of iterations to the normalized noise needs to be sufficiently small, capped by the specified privacy level~$\overline{\epsilon}_i$ and the agent's maximum capacity. It also shows the effect of the sensitivity on this trade-off.
The latter equation shows that with increasing number of agents~$M$ injecting noise, we need more iterations to achieve the same level of suboptimality. 
Similarly, if the maximum capacity~$u_{\text{max}}$ of the agents increases or the maximum price~$\pi_{\text{max}}$ decreases, we require more iterations or lower noise variance to maintain the same level of suboptimality.





\section{Sharing or Pricing the Privacy Budget?}
Theorem \ref{th:convergence rate of inexact AMA} provides a relationship between sub-optimality and the cumulative variance of the Laplacian noise inserted by all subsystems. 
In real scenarios, a system designer or operator may specify a desired level of (sub-)optimality $ \mathcal S_K$ achieved after $K$ iterations, that is $\mathbb{E}[D(\lambda^{K}) - D(\lambda^{\star})] \le \mathcal S_K$. 
Rewriting Equation~\eqref{eq: upper-bound of the distributed optimization algorithm}, we can compute a bound on the amount of cumulative Laplacian noise allowed at run-time
\begin{equation}
\sum^{M}_{i=1} \sigma^{2}_i \le  \frac{1}{4} \rho^{2}_{\phi} K \mathcal S_K - \sum^{M}_{i=1} G^{2}_i \triangleq \Sigma_{\text{budget}}  \ .
\end{equation}
Hence, once $\mathcal S_K$ is specified, a \emph{cumulative privacy budget} $\Sigma_{\text{budget}}$ is set.
Remember that each individual agent may define a different distance metric, or different weights, which may lead to different sensitivities and hence different noise levels required by the various agents.
Since the privacy budget is limited, a fair and transparent allocation procedure is required to divide the allowable noise over all agents.
Here, we propose two approaches for going about such allocation.

The first approach, would entail a proportional division of the pie. This could be done in two ways. The first way involves splitting the noise variance budget by the number of agents:~$\sigma_i^2 = \frac{1}{M} \Sigma_{\text{budget}} \,, \forall i = 1,\hdots,M$. In this case, the allowed noise variance $\sigma_i^2$ determines the maximum level of local differential privacy~$\epsilon_i$ that can be maintained, given a set sensitivity~$\Theta_i$ (or vice versa), as outlined in~\eqref{eq:privacy_level}), which will differ from agent to agent. The second way involves setting all local differential privacy levels~$\epsilon_i$ equal, and, given set sensitivities~$\Theta_i$, splitting the noise among all agents. This is equivalent to equating Equation \eqref{eq:privacy_level} for all agents $i = 1, \hdots,M$, leading to $M$ equations for~$\sigma_1,\hdots,\sigma_M$:
\begin{equation}
 \frac{\Theta_1}{\sigma_1} =  \frac{\Theta_2}{\sigma_2} = \cdots =
 \frac{\Theta_M}{\sigma_M}, \quad \sum^{M}_{i=1} \sigma^{2}_i = \Sigma_{\text{budget}} \,,
\end{equation}
where we assume that the noise variances are constant for all time steps~$k = 1,\hdots,K$.

The second approach we anticipate is a pricing scheme, in which the value of privacy is left to a market or negotiation.
In the context of d-OPF, it is natural to assume that different DER owners with varying privacy levels will have varying degrees of willingness to pay or incur a deduction on their revenue for preserving privacy of their local parameters.
Here, we propose two scenarios to perform allocation via the so-called \emph{Kelly mechanism}~\cite{kelly_rate_1998}.
We assume a one-directional bid $w_i$ done by all agents after seeing a price $\pi^{\sigma}_i$ given by the network operator to the agent. With the bid, the operator constructs a surrogate utility function~$w_i \log \sigma_i^2$. The operator then determines the allocation and payment by maximizing the sum of surrogate utility functions as:
\begin{equation}
\sigma^2_{\text{Kelly}}(w) = \arg \displaystyle \max_{ \Sigma_{i=1}^M \sigma_i^2 \le \Sigma_{\text{budget}}} \sum_{i=1}^M w_i \log \sigma_i^2 \, .
\end{equation}
The authors of~\cite{yang_vcg-kelly_2007} show that this mechanism works neatly if one assumes that buyers are \emph{price-takers}. 
In the case buyers are strategic, realizing that the eventual price is influenced by all bids, i.e. $\pi^{\sigma}_i(w)$, the Kelly mechanism may not yield an efficient Nash equilibrium.
To account for such behavior, the same authors combined the Kelly mechanism with the celebrated Vickrey-Clark-Groves (VCG) mechanism~\cite{vickrey_counterspeculation_1961,clarke_multipart_1971,groves_incentives_1973}, proposing the VCG-Kelly mechanism. 
The VCG mechanism is lauded because it incentivizes all participating buyers to report truthfully about the value they assign to the goods sold. It has an intuitive payment scheme, which says that each player $i$ should pay the difference in utility of the other players $j \neq i$ between the scenarios that player $i$ does and does not participate, which leads to each player caring about both its own utility and that of others, leading to truthful reporting.
The VCG-Kelly mechanism extends the VCG mechanism to problems with divisible goods~\cite{yang_vcg-kelly_2007}. It considers a more general class of surrogate utility functions~$V_i(w_i,\sigma^2_i) = w_i f_i(\sigma^2_i)$, with $f_i$'s strictly increasing, strictly concave, and twice differentiable. 
The allocation rule is then similar in form:
\begin{equation}
\sigma^2_{\text{VCGK}}(w) = \arg \displaystyle \max_{ \sum_{i=1}^M \sigma_i^2 \le \Sigma_{\text{budget}}} \sum_{i=1}^n V_i(w_i,\sigma^2_i) \, ,
\end{equation}
and the payment scheme reads:
\begin{equation}
\begin{array}{rcl}
m_{i,\text{VCGK}}(w) &=& \displaystyle \left( \max_{\sum_{i=1}^M \sigma_i^2 \le \Sigma_{\text{budget}}, \sigma^2_i=0} \displaystyle \Sigma_{j=1,j \neq i}^M V_j(w_j,\sigma^2_j) \right) \\
&&- \displaystyle \sum_{j=1,j \neq i}^M V_i(w_i,\sigma^2_{i,\text{VCGK}})\, .
\end{array}
\end{equation}
The actual form of the utility functions in the context of our application is left as an open problem, which requires further investigation. Readers interested in further details and examples of the VCG-Kelly mechanism are directed to~\cite{yang_vcg-kelly_2007}.

We propose these allocation schemes to trigger a discussion about what a fair and workable division of the allowable noise variance in our privacy-preserving looks like. 
This may vary based on the application. 
Given that privacy in engineered systems is a value of increasing importance, it may be wise to consider the impact of pricing privacy for participants with varying socioeconomic background. In situations where certain agents have significantly less resources but privacy is equally important, a proportional scheme may be the more ethical approach to take.
If pricing is still used and it is anticipated that some participants outbid others, one may think about allowing each participant to have a minimum amount of privacy, which can be translated into an extra linear constraint in the Kelly mechanisms.

\section{Conclusions and Future Work}

In this paper, we developed local $\epsilon$-differential privacy for distributed optimization, building on recent advances the inexact alternating minimization algorithm (IAMA).
By exploiting the IAMA's convergence properties under the existence of errors in communication and computation, we are able to add noise to agent-to-agent communication in a way that preserves privacy in the specifications of user objectives and constraints while still guaranteeing convergence.
The method extends current approaches for differential privacy in distributed optimization by allowing privacy for both objectives and constraints and customization of privacy specifications for individual agents.
We analyzed the trade-offs between privacy and suboptimality for various levels of noise and number of iterations, and gave a method to determine feasible values of noise variance and number of iterations given specifications on privacy and suboptimality. 
We propose different alternatives to allocate the allowable noise variance across participating agents, either via proportional sharing or market mechanisms that incentivize for truthful reporting and allow efficient Nash equilibria, which can be used to implement the algorithm in fair and efficient ways.

\bibliographystyle{abbrv}
\bibliography{references}


\section{Supplementary Material}

\noindent Supplementary material, covering proofs for theorem and lemmas, and numerical results for an application to distributed model predictive control.

\subsection{Proof of Theorem \ref{thm:privacy}}
\label{append:proof_privacy}
\begin{IEEEproof}
To show that the proposed mechanism $\mathcal A$  in Algorithm \ref{al: Distributed optimization with errors} has the promised privacy guarantee, we first show that $\mathcal A_i$ is locally $\epsilon_i$-differentially private for node $i$. To this end, we study the quantity of interest
\begin{align}
\frac{\pp{\mathcal A_i^K(\mathcal P_1,\ldots,\mathcal P_i,\ldots, \mathcal P_M)\in\mathcal Z}}{ \pp{\mathcal A_i^K(\mathcal P_1,\ldots,\mathcal P'_i,\ldots, \mathcal P_M)\in\mathcal Z}}
\end{align}
of node $i$. In the proof we use random variables $Z_i^k, V_i^k$, for $k=1,\ldots, K$, to denote the output of the local mechanism $\mathcal A_i$ (c.f. Eq. (\ref{eq:Ai_output})) with input $(\mathcal P_1,\ldots,\mathcal P_i,\ldots,\mathcal P_M)$, and use $Z_i^{'k},V_i^{'k}$, for $k=1,\ldots, K$, to denote the output of the local mechanism $\mathcal A_i$ with input $(\mathcal P_1,\ldots,\mathcal P_i',\ldots,\mathcal P_M)$ where the two local problems $\mathcal P_i$ and $\mathcal P_i'$ have distance $d(\mathcal P_i,\mathcal P_i')\leq 1$. With these notations, we can rewrite the above expressions as
\begin{align*}
\frac{\pp{ Z_i^1= z_i^1, V_i^1=v_i^1, \ldots, Z_i^K= z_i^K, V_i^K=v_i^K }}{\pp{ Z_i^{'1}= z_i^1, V_i^{'1}=v_i^1, \ldots, Z_i^{'K}= z_i^K, V_i^{'K}=v_i^K  }}
\end{align*}
for some realizations $z_i^k, v_i^k, k=1,\ldots K$.

It is straightforward to show that if $Z_i^k=Z_i^{'k}, V_i^k=V_i^{'k}$ for all $k=1,\ldots, K$, then we have $Z_j^k=Z_j^{'k}, V_j^k=V_j^{'k}$ for all $j\neq i$ and for all $k=1,\ldots, K$. For the simplicity of notation, we use $Z^k$ and $V^k$ to denote the tuple $(Z_1^k,\ldots, Z_M^k)$ and $(V_1^k,\ldots, V_M^k)$, respectively. Furthermore notice  that $V^k$ is a deterministic function of $Z^k$, as shown in the step 3 of Algorithm \ref{al: Distributed optimization with errors}. Hence we have
\begin{align*}
&\frac{\pp{ Z_i^1= z_i^1, V_i^1=v_i^1, \ldots, Z_i^K= z_i^K, V_i^K=v_i^K }}{\pp{ Z_i^{'1}= z_i^1, V_i^{'1}=v_i^1, \ldots, Z_i^{'K}= z_i^K, V_i^{'K}=v_i^K  }}\\
&=\frac{\pp{ Z^1= z^1, V^1=v^1, \ldots, Z^K= z^K, V^K=v^K }}{\pp{ Z^{'1}= z^1, V^{'1}=v^1, \ldots, Z^{'K}= z^K, V^{'K}=v^K  }}\\
&=\frac{\pp{ Z^1= z^1, \ldots, Z^K= z^K }}{\pp{ Z^{'1}= z^1 \ldots, Z^{'K}= z^K}}\\
&=\prod_{k=1}^K \frac{\pp{Z^k=z^k|Z^{t}=z^{t},t<k}}{\pp{Z^{'k}=z^k|Z^{'t}=z^{t},t<k}}
\end{align*}
where $z^k, v^k, k=1,\ldots, K$ are realizations of the random variables. Now we analyze the term
\begin{align*}
 \frac{\pp{Z^k=z^k|Z^{t}=z^{t},t<k}}{\pp{Z^{'k}=z^k|Z^{'t}=z^{t},t<k}}
\end{align*}
for each $k=1,\ldots, K$. 

According to Step 1 of Algorithm \ref{al: Distributed optimization with errors}, $z_i^k$ is determined by $\mu_i^{k-1}$ and $\delta_i^k$, in other words $Z^k - m^{k-1}-(Z_0,\ldots, Z^{k-1})$ forms a Markov chain where the random vector $m^k$ denotes the Lagrangian $(\mu_1^k,\ldots, \mu_M^k)$ in the algorithm. Consequently, it holds that
\begin{align*}
 \frac{\pp{Z^k=z^k|Z^{t}=z^{t},t<k}}{\pp{Z^{'k}=z^k|Z^{'t}=z^{t},t<k}}=\frac{\pp{Z^k=z^k | m^{k-1}=\mu^{k-1}}}{\pp{Z^{'k}=z^k | m^{'k-1}=\mu^{'k-1}}}
\end{align*}
where $m^k$ and $m^{'k}$ denote the random Lagrangians with inputs $(\mathcal P_1,\ldots,\mathcal P_i,\ldots,\mathcal P_M)$ and $(\mathcal P_1,\ldots,\mathcal P_i',\ldots,\mathcal P_M)$, respectively. Furthermore, it can be checked easily that $m^{k-1}$ and $m^{'k-1}$ are completely determined by $Z^t, Z^{'t}$ for $ t<k$.  
Since we have $Z^t=Z^{'t}=z^t$ for all $t<k$, we also have  $\mu^{k-1}=\mu^{'k-1}$. 

With the above observation, we can continue our derivation
\begin{align*}
& \frac{\pp{Z^k=z^k|Z^{t}=z^{t},t<k}}{\pp{Z^{'k}=z^k|Z^{'t}=z^{t},t<k}}=\frac{\pp{Z^k=z^k|m^{k-1}=\mu^{k-1}}}{\pp{Z^{'k}=z^k|m^{'k-1}=\mu^{'k-1}}}\\
 &=\frac{\pp{Z^k=z^k|m^{k-1}=\mu^{k-1}}}{\pp{Z^{'k}=z^k|m^{'k-1}=\mu^{k-1}}}\\
 &=\prod_{j=1}^M \frac{\pp{Z_j^k=z_j^k|m_j^{k-1}=\mu_j^{k-1}}}{\pp{Z_j^{'k}=z_j^k|m_j^{'k-1}=\mu_j^{k-1}}}\\
&=\prod_{j=1}^M\frac{\pp{g(\mathcal P_j,\mu_i^{k-1})+\delta_j^k=z_j^k}}{\pp{g(\mathcal P_j', \mu_i^{k-1})+\delta_j^k=z_j^k}}\\
&=\prod_{j=1}^M \frac{\exp \left({-\|z_j^k-g(\mathcal P_j, \mu_j^{k-1})\|/\sigma_j^k}\right)}{\exp \left(-\|z_j^k-g(\mathcal P_j', \mu_j^{k-1})\|/\sigma_j^k\right)}\\
&\leq \prod_{j=1}^M \exp\left(  \|g(\mathcal P_j, \mu_j^{k-1})-g(\mathcal P_j',\mu_j^{k-1}))\|/\sigma_j^k\right)\\
&=\exp\left(  \|g(\mathcal P_i,\mu_i^{k-1})-g(\mathcal P_i',\mu_i^{k-1}))\|/\sigma_i^k\right)\\
&\leq \exp (\Theta_i/\sigma_i^k)
\end{align*}
with $\Theta_i$ defined in (\ref{eq:sensitivity}). Finally, we obtain
\begin{align*}
&\frac{\pp{\mathcal A_i^K(\mathcal P_1,\ldots,\mathcal P_i,\ldots, \mathcal P_M)\in\mathcal Z}}{ \pp{\mathcal A_i^K(\mathcal P_1,\ldots,\mathcal P'_i,\ldots, \mathcal P_M)\in\mathcal Z}}\\
&=\prod_{k=1}^K \frac{\pp{Z^k=z^k|Z^{t}=z^{t},t<k}}{\pp{Z^{'k}=z^k|Z^{'t}=z^{t},t<k}}\\
&\leq \prod_{k=1}^K \exp ( \Theta_i/\sigma_i^k)\\
&=\exp \left(\sum_{k=1}^K \Theta_i /\sigma_i^k \right)
\end{align*}
This proves the privacy guarantee in (\ref{eq:privacy_level}).

\end{IEEEproof}

\subsection{Proofs of Lemmas}

\begin{IEEEproof}[Proof of Lemma \ref{lemma:independent_mu}]
We show that  for the quadratic problem, the optimization problem in (\ref{eq:sensitivity_quadratic}) is in fact independent of the value  of $\mu$.  To see this, we define $\tilde h_i:=h_i-\mu$ and $\tilde h_i':=h_i'-\mu$. The optimization problem (\ref{eq:sensitivity_quadratic}) can be rewritten as
\begin{align}
\text{maximize }&\Theta_i:= \|z^*-z'^{*}\| \label{pr:sensitivity}\\
\text{s. t. } &z^*\in \text{KKT}(H_i,\tilde h_i,C_i,c_i,w_i, 0) \nonumber\\
&z'^{*}\in\text{KKT}(H_i',\tilde h_i',C_i',c_i',w_i', 0) \nonumber\\
&a_1\|H_i-H_i'\|+a_2\|\tilde h_i-\tilde h_i'\|\nonumber\\
&+a_3\|C_i-C_i'\|+a_4\|c_i-c_i'\| \leq 1\nonumber\\
&\text{with variables } z^*,z'^{*}, H_i, H_i', \tilde h_i, \tilde h_i',\nonumber \\
&C_i, C_i', c_i, c_i', w_i, w_i' \nonumber
\end{align}
where we used the fact $\tilde h_i-\tilde h_i' = h_i-h_i'$. 
It can be seen that the above problem formulation is independent of $\mu$, which proves the claim.
\end{IEEEproof}

\begin{IEEEproof}[Proof of Lemma \ref{lemma:sensitivity_change_quadratic}]
Using Lemma \ref{lemma:independent_mu}, we consider the
  case with  $g(P_i, \mu=0)=\operatorname{argmin}_{z\in \mathbb C} \frac{1}{2}z^{T}Hz + h^{T}z$ (we omit the index $i$ for simplicity). Define 
\begin{align*}
z^* := \operatorname{argmin}_{z\in \mathbb{C}} \;\frac{1}{2}z^{T}Hz + h^{T}z
\end{align*}
and 
\begin{align*}
z'^* := \operatorname{argmin}_{z\in \mathbb{C}} \;\frac{1}{2}z^{T}H'z + h^{T}z
\end{align*}
where $H, H'$ satisfies $\|H-H'\|_2\leq 1$.

The optimality of $z^*$ and feasibility of $z'^*$ implies 
\begin{align*}
\langle Hz^*+h, z'^*-z^*\rangle\geq 0,
\end{align*}
and the optimality of $z'^*$  and the feasibility of $z^*$ implies that
\begin{align*}
\langle H'z'^{*} +h, z^*-z'^*\rangle \geq  0.
\end{align*}

Manipulations of the above two inequalities yield
\begin{align*}
\langle H'z'^*-Hz^*, z^*-z'^*\rangle\geq 0
\end{align*}
which can be rewritten as
\begin{align*}
\langle H'\Delta + (H-H')z^*, \Delta\rangle\leq 0
\end{align*}
where $\Delta:=z^*-z'^*$. This implies that
\begin{align*}
\Delta^TH'\Delta\leq \Delta^T(H'-H)z^*
\end{align*}
Using the fact that $\Delta^TH'\Delta\geq \lambda_{min}(H')\|\Delta\|^2$ and $ \Delta^T(H'-H)z^*\leq \|\Delta\|\|z^*\||\lambda_{max}(H'-H)|$, we conclude that 
\begin{align*}
\|z^*-z'^{*}\|&\leq \frac{|\lambda_{max}(H^{\prime}-H)|}{\lambda_{\min}(H^{\prime})} \cdot \|z^*\|\\
&\leq \frac{1}{\lambda_{min}^{(i)}}G_i
\end{align*}
where the last inequality uses the assumption that $|\lambda_{max}(H^{\prime}-H)|=\|H-H'\|_2\leq 1$, $\|z^*\|\geq G_i$, and the minimum eigenvalue of $H'$ is lower bounded by $\lambda_{min}^{(i)}$
\end{IEEEproof}

\begin{IEEEproof}[Proof of Lemma \ref{lemma:sensitivity_change_linear}]
Using Lemma \ref{lemma:independent_mu}, we consider the
  case with  $g(P_i, \lambda=0)=\operatorname{argmin}_{z\in \mathbb C} \frac{1}{2}z^{T}Hz + h^{T}z$ (we omit the index $i$ for simplicity). Define 
\begin{align*}
z^* := \operatorname{argmin}_{z\in \mathbb{C}} \;\frac{1}{2}z^{T}Hz + h^{T}z
\end{align*}
and 
\begin{align*}
z'^* := \operatorname{argmin}_{z\in \mathbb{C}} \;\frac{1}{2}z^{T}Hz + h'^{T}z
\end{align*}
where $h, h'$ satisfies $\|h-h'\|_2\leq 1$.

Since $H\succ 0$, we define $H = D\cdot D^{T}$ with $D$ invertible and rewrite the above expression as
\begin{align*}
z^*& =  \operatorname{argmin}_{z\in \mathbb{C}} \;\frac{1}{2}z^{T}Hz + h^{T}z\\
& = \operatorname{argmin}_{z\in \mathbb{C}} \;\frac{1}{2} \|D^{T}z + D^{-1}h\|^{2} 
\end{align*}
Let $v = D^{T}z$. The optimization problem above becomes 
\begin{align*}
v^*
& =  \operatorname{argmin}_{v\in \bar{\mathbb{C}}} \frac{1}{2} \|v + D^{-1}h\|^{2} \enspace ,
\end{align*}
where we define $\bar{\mathbb{C}}:=\{D^T z\mid z\in \mathbb{C}\}$.
Similarly we have
\begin{align*}
v'^*
& =  \operatorname{argmin}_{v\in \bar{\mathbb{C}}} \frac{1}{2} \|v + D^{-1}h'\|^{2} \enspace.
\end{align*}

Hence $v^*$ (res. $v'^*$) can be seen as the projection of the point $-D^{-1}h$ (res. $-D^{0-1}h'$)  onto the set $\bar{\mathbb{C}}$.  Since $\mathbb{C}$ is convex, we know that $\bar{\mathbb{C}}$ is also convex, and it holds that 
\begin{align*}
\|v^*-v'^*\| \leq \|D^{-1} (h-h')\| \enspace .
\end{align*}
Using $z = (D^T){^{-1}}v$, we obtain
\begin{align*}
\|z^*-z'^{*}\|_2 &\leq \|D^{-1}\|_2\cdot\|D^{-1}\cdot (h-h')\|_2 \\
& \leq \|D^{-1}\|_2^{2}\cdot\|h-h'\|_2 \\
&\leq \frac{1}{\lambda_{min}(H)}\cdot \|h-h'\|_2\\
&\leq \frac{1}{\lambda_{min}^{(i)}}
\end{align*}
where the last inequality uses the assumption that $\|h-h'\|\leq 1$ and  minimum eigenvalue of $H$ is lower bounded by $\lambda_{min}^{(i)}$.
\end{IEEEproof}

\begin{IEEEproof}[Proof Sketch for Lemma \ref{le:the equivalence between distributed optimization algorithm and stochastic PGM}]
This proof is an extension of the proof found for~\cite[Lemma~3.4]{pu_inexact_2016}. We first write Problem~\ref{pr:dist_opt} in the form of the splitting problem~\cite[Problem~3.1]{pu_inexact_2016}, by defining the two objectives as $f(z)=\sum^{M}_{i=1} f_{i}(z_i)$ subject to $z_i\in \mathbb{C}_{i}$ for all $i = 1,\cdots,M$ and $g=0$ with the optimization variables for the two objectives $z = [z^{T}_1,z^{T}_2,\cdots ,z^{T}_M]^{T}$ and $v$, respectively. The coupling matrices are set to $A=I$, $B=-E = -[E^{T}_1,E^{T}_2,\cdots ,E^{T}_M]^{T}$ and $c=0$. Under the assumption that the first objective $f(\mathbf{z})$ consists of a strongly convex function on $z$ and convex constraints. The convex constraints can be considered as indicator functions, which are convex functions. Due to the fact that the sum of a strongly convex and a convex function is strongly convex, the objective $f(z)$ is a strongly convex function. The second objective $g$ is a convex function. Then, by using the results in \cite[Lemma~3.4]{pu_inexact_2014}, we can prove that Algorithm~\ref{al: Distributed optimization with errors} is equivalent to Algorithm~\ref{al:stochastic PGM}, executed on Problem \ref{pr:dual problem of distributed optimization problem}, the dual problem of Problem~\ref{pr:dist_opt}.
\end{IEEEproof}

\section{Application: Distributed Model Predictive Control}
\label{sec:d-mpc}

\begin{figure*}[h!]
\centering
\subfloat[The sum of the inputs of all rooms, $\sum^{M}_{i = 1} u_{i}(t)$.]{%
  \includegraphics[width=8cm]{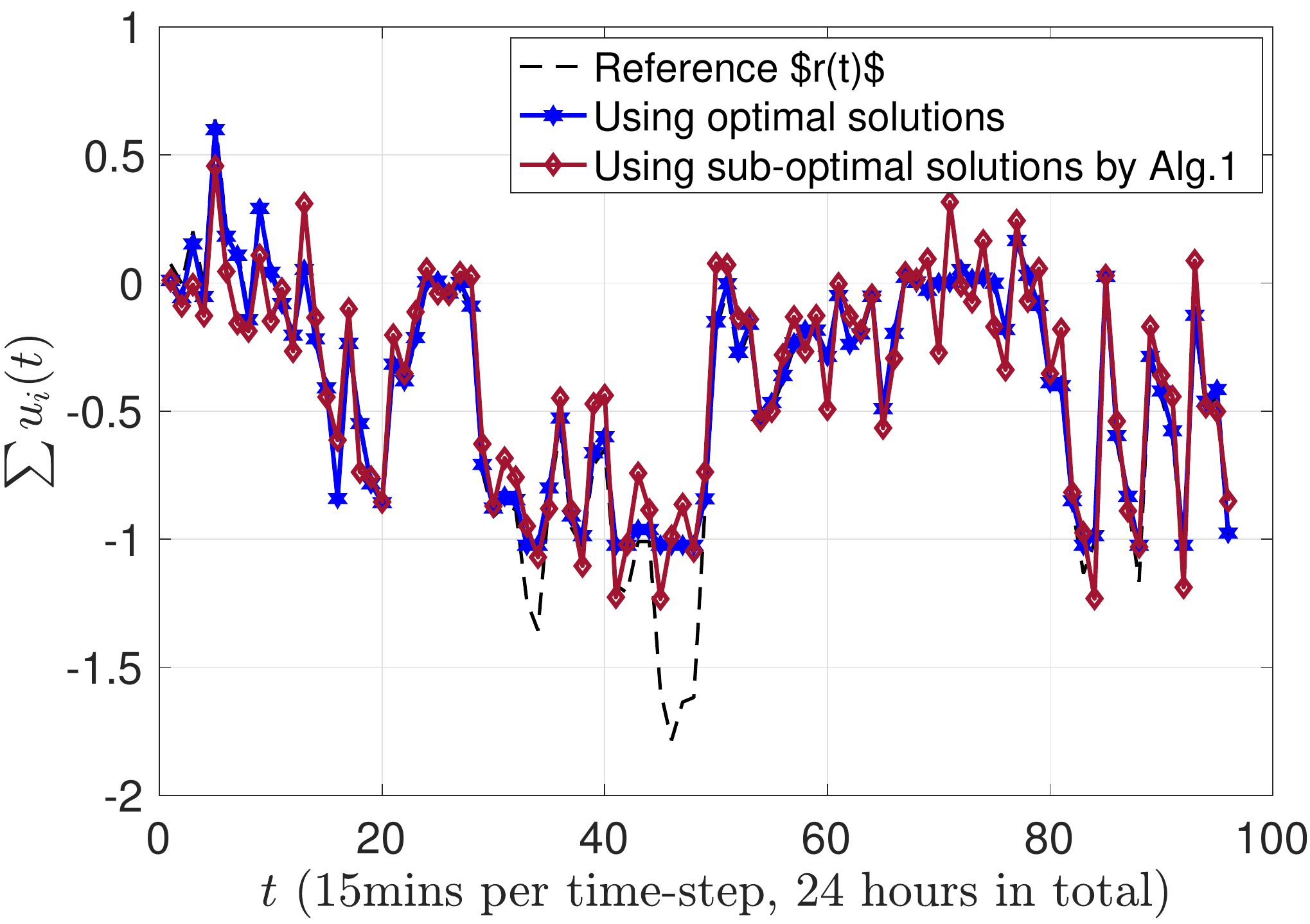}%
  \label{fig:close loop inputs tracking using DMPC}%
}\qquad
\subfloat[The temperature of room $1$.]{%
  \includegraphics[width=8cm]{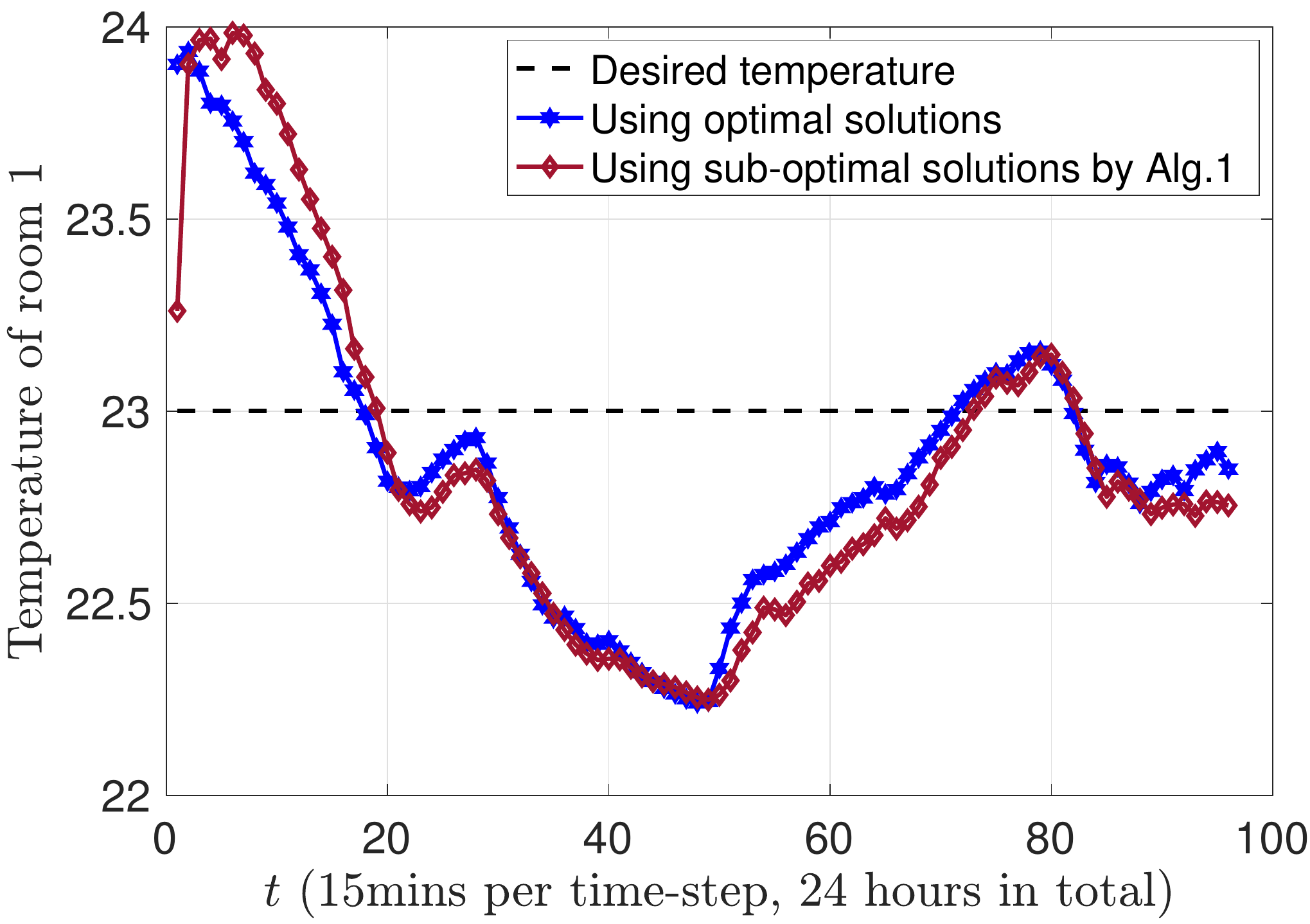}%
  \label{fig:close loop temperature of room 1 using DMPC}%
}
\caption{Simulation results}
\end{figure*}

In this section, we demonstrate the integration of customized local differential privacy in distributed model predictive control. We solve a problem that considers regulating the temperature in a set of rooms in a building, while making sure the cumulative control signal tracks a reference value set by a grid operator in a demand response scheme. 
In this example, we consider a linearized discrete-time dynamical model for each sub-system (room):
\begin{equation}
x_{i}(t+1) = A_{i}x_{i}(t)+B_{i}u_{i}(t)\enspace .
\end{equation}
Each sub-system has two states $x_i(t)\in\mathbb{R}^{2}$ and one input $u_i(t)\in\mathbb{R}^{1}$. The first state represents the temperature of room $i$ and the second state is a virtual state that is the result of linearization around an intended operating point.
The input $u_i$ represents the power input of the heating equipment in room $i$. The control goal of this example is keep of the temperature of each room at a desired level of $23^\circ$C and let the power consumption of all rooms to track a signal given by a higher level planner, denoted by $r(t)$. The dynamical matrices $A_i$ and $B_i$ are assumed to be controllable. We formulate the following problem:
\begin{problem}[Distributed Model Predictive Control]
\label{pr: MPC}
\begin{align*}
\min_{x,u}  & \quad \sum^{N}_{t=0} \sum^{M}_{i=1}  x_{i}(t)^{T} Q_i x_i(t) + u_{i}(t)^{T}  R_i  u_i(t) \\
&\quad\quad+ \alpha_{tr} (\sum^{M}_{i=1} u_i(t) - r(t))^{2} \\
s.t. \quad &x_{i}(t+1) = A_{i}x_{i}(t)+B_{i}u_{i}(t)\\
& u_i(t) \in \mathbb{U}_i, \; x_i(0) = \bar{x}_i, \quad i=1,2,\cdots ,M\enspace.
\end{align*}
\end{problem}
\noindent Here, the weight matrices are set as $Q_i = \mbox{diag}[1 0]\succeq 0$ and $R_i = \mbox{diag}[1 1]\succ 0$, and the tracking parameter $\alpha_{tr} = 10$. 
The horizon for the MPC problem is set to $N = 13$. 
The state and input sequences along the horizon of agent $i$ are denoted by $x_{i}=[x^{T}_{i}(0),x^{T}_{i}(1),\cdots,x^{T}_{i}(N)]^{T}$ and $u_{i}=[u^{T}_{i}(0),u^{T}_{i}(1),\cdots,u^{T}_{i}(N)]^{T}$. 
The input constraint $\mathbb{U}_i$ for sub-system $i$ is set to $\mathbb{U}_i = \{u_i| -0.2047 \leq u_i(t)\leq 0.7409\}$. This constraint guarantees that the linearized model is a reasonable approximation of the original non-linear model. 
The linearization allows us to formulate a distributed QP of the form in Problem~\ref{pr:dist_opt} with local variables $z_0 =  [u^{T}_{1},u^{T}_{2},\cdots,u^{T}_{M}]^{T}$, and $z_i = u_i$ for $i = 1, \ldots, M$. The global variable $v= [u^{T}_{1},u^{T}_{2},\cdots,u^{T}_{M}]^{T}$. The $i$-th component of $v$ is equal to $[v]_i = [u_{i}]$. 
Matrix $H_i$ is dense and positive definite, since the dynamical matrices $A_i$ and $B_i$ are controllable and the matrix $R_i$ is chosen to be positive definite. 
The local input constraints are $\mathbb{C}_i = \mathbb{U}_i$ for $i = 1, \ldots, M$, and $\mathbb{C}_0 = \mathbb{R}^M$. 

We consider the tracking signal $r(t)$ to be a confidential signal designed by the utility operator, containing sensitive market information. 
Therefore, in the resulted distributed QP, Problem~\ref{pr:dist_quad_opt}, the vector $h_i$ in the first objective is the private parameter to protect. 
By using the sampling rule of Lemma~\ref{le:sampling_rule} and the analytical upper bound (discussed in supplementary material), we can compute an upper-bound of the sensitivity of the distributed QP problem as $\Theta_i \leq 0.9756$. 
We test Algorithm~\ref{al: Distributed optimization with errors} for solving this distributed QP. 
At each time-step $t$, we utilize a warm-starting strategy to solve the MPC problems by taking the solution of the previous step to be the initial solution for the current step, and we set the number of iterations and the noise variance to be $K = 10$ and $\sigma = 0.1$. 
We compare the performance of the closed-loop system of the sub-optimal solutions given by Algorithm~\ref{al: Distributed optimization with errors} to the optimal solutions, over a period of $24$ hours ($96$ steps, $15$ mins per step). Fig.~\ref{fig:close loop inputs tracking using DMPC} shows the tracking performance of the sum of the inputs $u_i(t)$ of all rooms. Fig.~\ref{fig:close loop temperature of room 1 using DMPC} illustrates how the temperature of room $1$ varies. 
\ifCLASSOPTIONcaptionsoff
  \newpage
\fi

\end{document}